\documentclass[11pt]{article}
\usepackage{amsmath, amssymb, amsthm}

\theoremstyle{plain}

\numberwithin{equation}{section}

\title{On higher-power moments of $\Delta(x)$(II) }
\date{}

\begin{document}

\maketitle

\vspace{-1.5cm}

\begin{center}

\medskip

Wenguang  Zhai

\medskip

School of Mathematical Sciences \\
Shandong Normal University \\
Jinan 250014, Shandong \\
 P. R. China\\
zhaiwg@hotmail.com
\end{center}

\begin{center}{
Acta Arith. Vol.{\bf 114}(2004), 35-54}\end{center}

\begin{abstract}
 Let $\Delta(x)$ be the error term of the Dirichlet
divisor problem.   The  asymptotic formula of the integral
$\int_1^T\Delta^k(x)dx$ is established for any integer $3\leq k\leq
9$ by an unified method.  Similar results  are also established for
some other well-known error terms in the analytic number theory .
\end{abstract}


\section{\bf  Introduction and  main results }

\subsection{\bf Notations}\

Throughout this paper, let $d(n)$ denote the Dirichlet divisor
function,
 $r(n)$ denote the number of ways $n$ can be written as $n=x^2+y^2$ for
$x,y\in {\Bbb Z},$ and $a(n)$ denote the Fourier coefficients of a
holomorphic cusp form of weight $\kappa=2n\geq 12$ for the full
modular group , $\tilde a(n):=a(n)n^{-\kappa/2+1 /2}.$ For short, we
use $d, r, a, \tilde a$ denote these functions, respectively.
$\zeta(s)$ denotes the Riemann zeta-function.

Suppose $x>0, t>0.$ Define
\begin{eqnarray}
&&\Delta(x): =\sum_{n\leq x}d(n)-x\log x-
(2\gamma-1)x,\\
&&P(x): =\sum_{n\leq x}r(n)-\pi x,\\
&&A(x): =\sum_{n\leq x}a(n),\\
&&E(t):=\int_0^t|\zeta(\frac{1}{2}+iu)|^2du-t\log(t/2\pi)-(2\gamma-1)t.
\end{eqnarray}

Suppose $f:{\Bbb N}\rightarrow {\Bbb R}$ is any function such that
$f(n)\ll n^\varepsilon,$ $k\geq 2$ is a fixed integer. Define
\begin{equation}
s_{k;l}(f):=\sum_{\sqrt{n_1}+\cdots
+\sqrt{n_l}=\sqrt{n_{l+1}}+\cdots +\sqrt{n_k }} \frac{f(n_1)\cdots
f(n_k)}{(n_1\cdots n_k)^{3/4}}\hspace{3mm}(1\leq l<k),
\end{equation}
\begin{equation}
B_k(f):=\sum_{l=1}^{k-1}{k-1\choose
l}s_{k;l}(f)\cos\frac{\pi(k-2l)}{4}.
\end{equation}
We shall use $s_{k;l}(f)$  to denote both of the series (1.5) and
its value. We will prove  the convergence of $s_{k;l}(f)$ in Section
3.

Suppose $A_0>2$ is a real number, define
\begin{eqnarray*}
&&K_0: =\min\{n\in {\Bbb N}:n\geq A_0, 2|n\},\\
&&b(k):=2^{k-2}+\frac{k-6}{4},\\
&&\sigma(k,A_0):= \left\{\begin{array}{ll}
1/4,&\mbox{if $k-1<A_0/2$ ,}\\
\frac{A_0-k}{2(A_0-2)},& \mbox{if $ A_0/2+1\leq k<A_0,$}
\end{array}\right.\\
&&\delta_1(k,A_0): =\sigma(k,A_0)/2b(K_0),\\
&&\delta_2(k,A_0): =\frac{\sigma(k,A_0)}{2b(k)+2\sigma(k,A_0)}.
\end{eqnarray*}

 ${\Bbb N}$ denotes the set of all natural numbers.
$\varepsilon$ always denotes a sufficiently small positive constant
which may be different at different places.
 We will use the inequality  $d(n)\ll n^\varepsilon $ freely.
$SC(\Sigma)$ denotes the summation condition of the sum $\Sigma.$
$\mu(n)$ is the M\"obius function.

\subsection{\bf Introduction}\

In this paper we shall study the higher-power moments of $\Delta(x),
$ $P(x),$  $ A(x)$ and $E(t).$

We begin with the Dirichlet divisor problem. Dirichlet first proved
that
 $\Delta(x)=O(x^{1/2}).$
 The exponent $1/2$ was improved by many authors.
The latest result reads
\begin{equation}
\Delta(x)\ll x^{23/73}(\log x)^{315/146},
\end{equation}
which can be found in  Huxley[6]. It is conjectured that
\begin{equation}
\Delta(x)=O(x^{1/4+\varepsilon}),
\end{equation}
which is supported by the classical mean-square  result
\begin{equation}
\int_1^T\Delta^2(x)dx=\frac{(\zeta(3/2))^4}{6\pi^2 \zeta(3)}T^{3/2}
+O(T\log^5 T)
\end{equation}
proved by Tong[17] and the upper bound estimate
\begin{equation}
\int_1^T|\Delta(x)|^{A_0}dx\ll T^{1+A_0/4+\varepsilon},
\end{equation}
where $A_0>2$ is a fixed real number. The estimate of type (1.10)
can be found in Ivi\'c[7, Thm. 13.9] with $A_0=35/4$ and
Heath-Brown[5] with $A_0=28/3.$
 On the other hand , Voronoi[19] proved that
\begin{equation}
\int_1^T\Delta(x)dx=T/4+O(T^{3/4}),
\end{equation}
which in conjunction with (1.9) shows that $\Delta(x)$ has a lot of
sign change s and cancelations between the positive and negative
portions.

Tsang[18] first  studied the third- and fourth-power moments of
$\Delta(x).$ He proved that ( with  notations in  Section 1.1)
\begin{eqnarray}
&&\int_1^T\Delta^3(x)dx=\frac{3s_{3;1}(d)}{28\pi^3}T^{7/4}
+O(T^{7/4-1/14+\varepsilon}),\\
&&\int_1^T\Delta^4(x)dx=\frac{3s_{4;2}(d)}{64\pi^4}T^2+O(T^{2-1/23+\varepsilon})
.
\end{eqnarray}

Heath-Brown[5] proved that
 for $k=3,4,5,6,7,8,9$ the limit
$$\lim_{T\rightarrow\infty}T^{-1-k/4}\int_1^T \Delta(x)^kdx$$ exists.

In [20] the author improved Tsang's method and  proved that
\begin{eqnarray}
&&\int_1^T\Delta^3(x)dx=\frac{3s_{3;1}(d)}{28\pi^3}T^{7/4}
+O(T^{3/2+\varepsilon}),\\
&&\int_1^T\Delta^4(x)dx=\frac{3s_{4;2}(d)}{64\pi^4}T^2+O(T^{2-2/41}),\\
&&\int_1^T\Delta^5(x)dx=\frac{5(2s_{5;2}(d)-s_{5;1}(d))}{288\pi^5}T^{9/4}
+O(T^{9/4-5/816}).
\end{eqnarray}

But the argument of [20] fails for $k\geq 6.$

\subsection{\bf New results on higher-power moments of
$\Delta(x)$}\

In this paper we shall use a different approach to study the
higher-power moment s of $\Delta(x).$ This leads to the asymptotic
formulas of the integral $\int_1^T\Delta^k(x)dx$ for $3\leq  k\leq
9.$   Furthermore, if the estimate (1. 8) is true, then our approach
can give the   asymptotic formulas of $\int_1^T\Delta^k(x)dx$  for
any $k\geq 10.$

{\bf Theorem 1.} Let $A_0>9$ be a real number such that (1.10) holds
, then for any integer $3\leq k<A_0,$
 we have the asymptotic formula
\begin{equation}
\int_1^T\Delta^k(x)dx=\frac{B_k(d)}{(1+k/4)2^{3k/2-1}\pi^k}T^{1+k/4}
+O(T^{1+k/4-\delta_1(k,A_0)+\varepsilon})
\end{equation}

{\bf Remark 1.1.} From Ivi\'c's argument[7, Thm 13.9],
 we know that the value of $A_0$ for which (1.10) holds
depends on the large-value estimate and the upper bound estimate of
 $\Delta(x).$
If we insert the estimate (1.7) into the argument of Ivi\'c, we get
that (1.10) holds with $A_0=184/19.$ Whence for $k\in
\{3,4,5,6,7,8,9\},$ we can get the asymptotic formula (1.17).
Moreover,  if the estimate $\Delta(x)\ll x^{5/16-\delta}$ holds for
some small $\delta>0,$ then the asymptotic formula (1.17) holds for
$k =10.$

{\bf Remark 1.2.} For $k\geq 10,$ Theorem 1 is only an conditional
result.
 However, it tells
 us that for any $k\geq 10,$ the main term in the asymptotic formula of
$\int_1^T\Delta^k(x)dx$(if it exists) must have the form stated in
(1.17).

{\bf Remark 1.3.} We can state the following three conjectures about
$\Delta(x). $

Conjecture 1: The estimate (1.8) is true.

Conjecture 2: The estimate (1.10) is true for any $A_0>2$.

Conjecture 3: For any fixed $k\geq 3,$ there exists a constant
$\delta_k>0$ such
 that
the asymptotic formula
$$\int_1^T\Delta^k(x)dx=\frac{B_k(d)}{(1+k/4)2^{3k/2-1}\pi^k}T^{1+k/4}
+O(T^{1+k/4-\delta_k+\varepsilon})$$ holds.

It is well-known that Conjecture 1 and Conjecture 2 are equivalent.
From Theorem 1 we know that actually the three Conjectures are
equivalent. It is easy  to deduce Conjecture 2 from Conjecture 3. To
deduce Conjecture 3 from Conjecture 2, we take $A_0=2(k-1)$ and
$\delta_k=\delta_1(k,2(k-1)).$

{\bf Remark 1.4.} From (1.11) we know that the integral
$\int_1^T\Delta(x)dx$ have many cancelations from the positive and
negative portions of $\Delta(x)$. However, from (1.12) Tsang[18]
observed that this is not so for $\int_1^T\Delta^3(x)dx.$ From
Theorem 1  we know this is also not so for
$\int_1^T\Delta^k(x)dx\hspace{1mm}(k=5,7,9)$ since numerical
computation tells $B_k(d)>0$ for $k=5,7,9.$ Maybe $B_k(d)>0$ holds
for any odd $k\geq 3.$

\bigskip

 The constant
$\delta_1(k,A_0)$ is small for $k$ small.  If we combine Ivi\'c's
argument in the proof of Theorem 1,  we can get  the following
Theorem 2 for $3\leq k\leq 9$. Note that  the results for $k=3,4$
are  weaker than those of [20]. Theorem 2 for $k=5$ improved (1.16).

{\bf Theorem 2.}  For $3\leq k\leq 9,$ the asymptotic formula (1.17)
holds with $\delta_1(k,A_0)$ replaced by
 $\delta_2(k,184/19).$

Especially for $k=5, 6, 7 , 8, 9,$ we have
\begin{eqnarray}
&&\int_1^T\Delta^5(x)dx=\frac{5(2s_{5;2}(d)-s_{5;1}(d))}{288\pi^5}T^{9/4}
+O(T^{9/4-1/64+\varepsilon}),\\
&&\int_1^T\Delta^6(x)dx=\frac{5s_{6;3}(d)-3s_{6;1}(d)}{320\pi^6}T^{5/2}
+O(T^{5/2-35/4742+\varepsilon}),\\
&&\int_1^T\Delta^7(x)dx=\frac{7(5s_{7;3}(d)-3s_{7;2}(d)-s_{7;1}(d))}{2816\pi^7}T
^{11/4}
+O(T^{11/4-17/6312+\varepsilon}),\\
&&\int_1^T\Delta^8(x)dx=\frac{7(5s_{8;4}(d)-4s_{8;2}(d))}{6144\pi^8}T^{3}
+O(T^{3-8/9433+\varepsilon}),\\
&&\int_1^T\Delta^9(x)dx=\frac{3(3s_{9;1}(d)-12s_{9;2}(d)-28s_{9;3}(d)
+42s_{9;4}(d))}{26624\pi^9}T^{13/4}
\end{eqnarray}
\begin{eqnarray*}+O(T^{13/4-13/75216+\varepsilon}).\end{eqnarray*}

\subsection{\bf Higher-power moments of $P(x), A(x)$ and $E(t)$}\

The method  of proving  Theorem 1 and Theorem 2 can also be applied
to study the  higher-power moments of $P(x), A(x)$ and $E(t).$

\bigskip

The conjectured bound of $P(x)$ is
\begin{equation}
P(x)=O(x^{1/4+\varepsilon}),
\end{equation}
which is supported by
\begin{equation}
\int_2^TP^2(x)dx=(\frac{1}{3\pi^2}\sum_{n=1}^\infty
r^2(n)n^{-3/2})T^{3/2} +O(T\log^2 T)
\end{equation}
proved by Katai[14]. Tsang[18] also studied the third- and the
fourth-power moments of $P(x).$ His results were improved in the
author[20]. An asymptotic formula for the fifth-power moment of
$P(x)$ was also obtained in [20], which is further improved by the
following Theorem 3($k=5$.)

{\bf Theorem 3.} Let $A_0>9$ be a real number such that the estimate
\begin{equation}
\int_1^T|P(x)|^{A_0}dx\ll T^{1+A_0/4+\varepsilon}
\end{equation}
 is true ,
then for any integer  $3\leq k<A_0,$
 the asymptotic formula
\begin{equation}
\int_1^TP^k(x)dx=\frac{(-1)^kB_k(r)}{(1+k/4)2^{k-1}\pi^k}T^{1+k/4}
+O(T^{1+k/4-\delta_1(k,A_0)+\varepsilon})
\end{equation}
holds.

 Especially   for $3\leq k\leq 9,$ the asymptotic formula (1.26) holds with
$\delta_1(k,A_0)$ replaced by
 $\delta_2(k,184/19).$

{\bf Remark 1.5.} Ivi\'c [7, Thm 13.12] proved that the estimate
(1.25) holds
 for $A_0=35/4.$ If we insert the estimate $P(x)=O(x^{23/73+\varepsilon})$
(see Huxley[6])
 into his argument,
we find that (1.25) holds for $A_0=184/19.$

\bigskip

It is well-known that $A(x)$ has no main term and
 $A(x)\ll x^{\kappa/2-1/6+\varepsilon}.$
From Deligne[4], we have $|\tilde a(n)|\leq d(n).$

The conjectured bound of $A(x)$ is $A(x)\ll
x^{\kappa/2-1/4+\varepsilon}.$ Ivi\'c[9] proved that
\begin{equation}
\int_1^T A^2(x)dx=B_2T^{\kappa+1/2}+O(T^\kappa\log^5 T),
\end{equation}
where
$$B_2=\frac{1}{4\kappa+2}\sum_{n=1}^\infty a^2(n)n^{-\kappa-1/2}.$$
Ivi\'c[9] also proved that the estimate
\begin{equation}
\int_1^T|A(x)|^{A_0}dx\ll T^{1+A_0(2\kappa-1)/4+\varepsilon}
\end{equation}
holds for $A_0=8.$ Cai [3] studied the third- and fourth-power
moments of $A(x).$ His results were improved in the author[20]. In
[20] an asymptotic formula for the fifth-pow er moment of $A(x)$ was
also obtained, which is further improved by the case $k=5$ of the
following Theorem 4.

{\bf Theorem 4.} Let $A_0\geq 8$ be a real number such that (1.28)
 is true ,
then for any $3\leq k <A_0,$
 the asymptotic formula
\begin{equation}
\int_1^TA^k(x)dx=\frac{B_k(\tilde
a)}{(1+\frac{k(2\kappa-1)}{4})2^{3k/2-1}\pi^k}
T^{1+\frac{k(2\kappa-1)}{4}}
+O(T^{1+\frac{k(2\kappa-1)}{4}-\delta_1(k,A_0)+\varepsilon})
\end{equation}
holds.

Especially for $3\leq k\leq 7,$  the asymptotic formula (1.29) holds
with $\delta_1(k,A_0)$ replaced by $\delta_2(k,8).$

\bigskip

Many results for $E(t)$ parallel to those for $\Delta(x)$ have been
obtained; se e Ivi\'c[8] for a survey. The conjectured bound for
$E(t)$ is $E(t)\ll t^{1/4+\varepsilon},$ which is supported by
\begin{equation}
\int_2^TE^2(t)dt=\frac{2\zeta^4(3/2)}{3\zeta(3)\sqrt{2\pi}}T^{3/2}
+O(T\log^5 T),
\end{equation}
proved by Meurman[15]. It has been proved that (see Huxley[6])
\begin{equation}
E(t)\ll t^{72/227}(\log t)^{629/227}, \hspace{2mm}t>2.
\end{equation}
Ivi\'c[7, Thm. 15.7] proved that the estimate
\begin{equation}
\int_1^T|E(t)|^{A_0}dt\ll T^{1+A_0/4+\varepsilon}
\end{equation}
holds for $A_0=35/4.$ Inserting the estimate (1.31) into Ivi\'c's
argument, we find that (1.32) is true for $A_0=576/61.$

Tsang[18] studied the third- and fourth-power moment of $E(t)$ by
using the Atkinson's formula[1]. His results were further improved
by Ivi\'c[10]  in a different way. The author[20] obtained new
results on the third
 and the
fourth power moments  of $E(t).$ An asymptotic formula for the fifth
power moment of $E(t)$ was also obtained in [20], which is further
improved by the case $k=5$ of
 the
following Theorem 5..

{\bf Theorem 5.} Let $A_0>9$ be a real number such that the
estimates (1.10) and

(1.32)
 hold ,
then for any $3\leq k<A_0,$
 we have the asymptotic formula
\begin{equation}
\int_1^TE^k(t)dt=\frac{B_k(d)}{(1+k/4)2^{3k/4-1}\pi^{k/4}}
T^{1+k/4}+O(T^{1+k/4-\delta_1(k,A_0)+\varepsilon}).
\end{equation}

Especially for
 $3\leq k\leq 9,$ the asymptotic formula (1.33) holds with
$\delta_1(k,A_0)$ replaced by $\delta_2(k,576/61).$

{\bf Acknowledgement.} The author deeply thanks the referee for his
 warm and valuable comments.
The author also thanks Prof. Isao Wakabayashi for helpful
discussions , through which the proof of Lemma 3.1 was further
improved.

\section{\bf Some Preliminary Lemmas}

We need the following Lemmas.

{\bf Lemma 2.1.} The square-roots of squarefree numbers are linearly
independent over the integers.

\begin{proof} This is a classical result of Besicovitch[2]. \end{proof}

{\bf Lemma 2.2.} Suppose $k\geq 3,$  $(i_1,\cdots, i_{k-1})\in
\{0,1\}^{k-1}$ such that
$$\sqrt{n_1}+(-1)^{i_1}\sqrt{n_2}+(-1)^{i_2}\sqrt{n_3}+\cdots
+(-1)^{i_{k-1}}\sqrt{n_k} \not= 0.$$ Then we have
$$
|\sqrt{n_1}+(-1)^{i_1}\sqrt{n_2}+(-1)^{i_2}\sqrt{n_3}+\cdots
+(-1)^{i_{k-1}}\sqrt{n_k}| \gg \max(n_1,\cdots,
n_k)^{-(2^{k-2}-2^{-1})}.$$

\begin{proof}
The cases $k=3, 4$ are  Lemma  1  and  Lemma 2 of Tsang[18],
respectively. The proof for the general case is the same as  the
proof of Lemma 1 of [18]. We note that Heath-Brown[5] stated a
similar result for $k$ even.
\end{proof}

{\bf Lemma 2.3.} Suppose $A,B\in {\Bbb R}, A\not= 0,$ then
$$\int_T^{2T}\cos(A\sqrt t+B)dt\ll  T^{1/2}|A|^{-1}.$$

{\bf Lemma 2.4.} Suppose $k\geq 3,$  $(i_1,\cdots, i_{k-1})\in
\{0,1\}^{k-1},$ $(i_1,\cdots, i_{k-1})\not= (0,\cdots,0),$ $1<N_1,
N_2,\cdots ,N_k, $ $0<\Delta\ll E^{1/2},$ $E=\max(N_1, N_2,\cdots
,N_k).$
 Let
$${\cal A}={\cal A}(N_1, N_2,\cdots ,N_k; i_1,\cdots, i_{k-1}; \Delta)$$
denote the number of solutions of the inequality
\begin{equation}
|\sqrt{n_1}+(-1)^{i_1}\sqrt{n_2}+(-1)^{i_2}\sqrt{n_3}+\cdots
+(-1)^{i_{k-1}}\sqrt{n_k}| <\Delta
\end{equation}
with $N_j<n_j\leq 2N_j, 1\leq j\leq k.$ Then
$${\cal A}\ll \Delta E^{-1/2}N_1N_2\cdots N_k+ E^{-1}N_1N_2\cdots N_k.$$

\begin{proof}
Without loss of generality, suppose $E=N_k.$ If $(n_1,\cdots,n_k)$
satisfies (2. 1), then
$$\sqrt{n_1}+(-1)^{i_1}\sqrt{n_2}+(-1)^{i_2}\sqrt{n_3}+\cdots
+(-1)^{i_{k-2}}\sqrt{n_{k-1}}=(-1)^{i_{k-1}+1}\sqrt{n_k}+\theta
\Delta$$ for some $|\theta|<1.$ Whence we get
$$\left(\sqrt{n_1}+(-1)^{i_1}\sqrt{n_2}+(-1)^{i_2}\sqrt{n_3}+\cdots
+(-1)^{i_{k-2}}\sqrt{n_{k-1}}\right)^{2}= n_k+O(\Delta N_k^{1/2}).$$
Hence for fixed $(n_1,\cdots,n_{k-1}),$ the number of $n_k$ is $\ll
1+\Delta N_k^{1/2}$ and thus
$${\cal A}\ll \Delta N_k^{1/2}N_1N_2\cdots N_{k-1}+ N_1N_2\cdots N_{k-1}.$$
\end{proof}

\section{\bf On the series $s_{k;l}(d)$}

In this section we shall study the series $s_{k;l}(d).$ Suppose
$y>1$ is a large parameter, and define
$$s_{k;l}(d;y): =
\sum_{\stackrel{\sqrt{n_1}+\cdots +\sqrt{n_l}=\sqrt{n_{l+1}}+\cdots
+\sqrt{n_k}} {n_1,\cdots,n_k\leq y}} \frac{d(n_1)\cdots
d(n_k)}{(n_1\cdots n_k)^{3/4}},1\leq l<k.$$ We shall prove the
following Lemma  3.1.

{\bf Lemma 3.1.} We have
$$|s_{k;l}(d)-s_{k;l}(d;y)|\ll y^{-1/2+\varepsilon}, 1\leq l<k.$$

{\bf Remark.} Lemma 3.1 is still true if  the divisor function $d$
is replaced b y any function $f: {\Bbb N}\rightarrow{\Bbb R}$ with
$f(n)\ll n^\varepsilon.$

\begin{proof}
We shall prove Lemma 3.1 by induction in $k.$ The case $k=2$ is
easy. The case $k=3$ is contained already in Page 70 of Tsang[18],
later we always  suppose $k\geq 4.$ Since $s_{k;l}(d)=s_{k;k-l}(d),$
we suppose $l\leq k/2.$

By  the symmetry, we get
\begin{equation}
|s_{k;l}(d)-s_{k;l}(d;y)|\ll \sum_{\stackrel{\sqrt{n_1}+\cdots
+\sqrt{n_l}=\sqrt{n_{l+1}}+\cdots +\sqrt{n_k}} {n_1>y}}
\frac{d(n_1)\cdots d(n_k)}{(n_1\cdots n_k)^{3/4}}
\end{equation}
$$\ll U_1(d;y)+U_2(d;y),$$
say,  where
\begin{eqnarray*}
&&U_1(d;y):=\sum_{j=l+1}^k\sum_{\stackrel{\sqrt{n_1}+\cdots
+\sqrt{n_l} =\sqrt{n_{l+1}}+\cdots +\sqrt{n_k}} {n_1=n_j>y}}
\frac{d(n_1)\cdots d(n_k)}{(n_1\cdots n_k)^{3/4}},\\
&&U_2(d;y):=\sum_{\stackrel{\sqrt{n_1}+\cdots +\sqrt{n_l}
=\sqrt{n_{l+1}}+\cdots +\sqrt{n_k}} {n_1>y, n_1\not= n_j,l+1\leq
j\leq k}} \frac{d(n_1)\cdots d(n_k)}{(n_1\cdots n_k)^{3/4}}.
\end{eqnarray*}

If $l=1,$  then obviously $U_1(d;y)=0$. If $l>1,$ then by induction
we get
\begin{equation}
U_1(d;y)\ll \sum_{n>y}\frac{d^2(n)}{n^{3/2}}s_{k-2;l-1}(d)\ll
y^{-1/2+\varepsilon}.
\end{equation}

Now we estimate $U_2(d;y).$ Let
$I=\{1,\cdots,l\},J=\{l+1,\cdots,k\}.$ Suppose $(n_1,\cdots, n_k)\in
{\Bbb N }^k$ such that
$$(*):\hspace{5mm}
\sqrt{n_1}+\cdots +\sqrt{n_l}=\sqrt{n_{l+1}}+\cdots +\sqrt{n_k},
n_1\not= n_j,l+1\leq j\leq k.$$ Then there exist two sets
$I_0\subset I,J_0\subset J$ which satisfy the following properties:

(1). $1\in I_0;$

(2). $\sum_{i\in I_0}\sqrt{n_i}=\sum_{j\in J_0}\sqrt{n_j}$;

(3). For any real subset $I_0^{'}\subset I_0,J_0^{'}\subset J_0,$ we
have
 $$\sum_{i\in I_0^{'}}\sqrt{n_i}\not= \sum_{j\in J_0^{'}}\sqrt{n_j}.$$
If $(I_0,J_0)=(I,J),$ then we say $(n_1,\cdots, n_k)$ is a primitive
$(k,l)-$poi nt. Let ${\Bbb N}_{k;l}$ denote the set of all points in
${\Bbb N}^k$ which satisfy (*) and ${\Bbb N}_{k;l}^{*}$ the set of
all primitive $(k,l)-$points, respectively.
 Let ${\cal G}_{k;l}$
denote the set of all possible pairs $(I_0,J_0)$ when
$(n_1,\cdots,n_k)$ runs through ${\Bbb N}_{k;l}.$ Note that if
$l=1,$ then ${\cal G}_{k;l}=\{(I,J)\}.$

Suppose $(I_0,J_0)\in {\cal G}_{k;l}.$ Let $l_1=\#I_0, l_2=l-l_1,$
$k_1=\#I_0+\#J_0,k_2=k-k_1.$ From (*), we know that $k_1\geq 3.$
Define
\begin{eqnarray*}
&&R_{1}^{(I_0,J_0)}(d;y):=\sum_{\stackrel{\sqrt {n_1}+\cdots +\sqrt
{n_{l_1}}= \sqrt {n_{l_1+1}}+\cdots +\sqrt {n_{k_1}}}
{n_1>y,(n_1,\cdots, n_{k_1})\in {\Bbb
N}_{k_1;l_1}^{*}}}\frac{d(n_1)\cdots d(n_{ k_1})} {(n_1\cdots
n_{k_1})^{3/4}}.
\end{eqnarray*}

If $(I_0,J_0)\not=(I,J),$ then $l_1<l, k_1<k$ and we define
\begin{eqnarray*}
&&R_{2}^{(I_0,J_0)}(d):=\sum_{\sqrt {m_1}+\cdots +\sqrt {m_{l_2}}=
\sqrt {m_{l_2+1}}+\cdots +\sqrt {m_{k_2}}}\frac{d(m_1)\cdots
d(m_{k_2})} {(m_1\cdots n_{k_2})^{3/4}}.
\end{eqnarray*}
By the induction assumption, $R_{2}^{(I_0,J_0)}(d)\ll 1.$

If$ (n_1,\cdots, n_{k_1})\in {\Bbb N}_{k_1;l_1}^{*},$ then by Lemma
2.1 we have
$$n_j=s_j^2h,s_1+\cdots+s_{l_1}=s_{l_1+1}+\cdots+s_{k_1}, \mu(h)\not= 0.$$
$n_1>y$ implies that there exists at least one $n_j(l_1+1\leq j\leq
j_1)$ such t hat $n_j\gg y.$ We suppose $n_{k_1}\gg y.$ So we have
\begin{eqnarray*}
R_{1}^{(I_0,J_0)}(d;y) &&\ll \sum_{h} \sum_{\stackrel{s_1+\cdots
+s_{l_1}= s_{l_1+1}+\cdots +s_{k_1}}{s_1^2h>y,s_{k_1}^2h\gg
y}}\frac{d(s_1^2h) \cdots d(s_{k_1}^2h)}
{h^{3k_1/4}(s_1\cdots s_{k_1})^{3/2}}\\
&&\ll \sum_{h} \sum_{\stackrel{s_1+\cdots +s_{l_1}= s_{l_1+1}+\cdots
+s_{k_1}}{s_1^2h>y,s_{k_1}^2h\gg y}}\frac{d^2(s_1) \cdots
d^2(s_{k_1})d^{k_1}(h)}
{h^{3k_1/4}(s_1\cdots s_{k_1})^{3/2}}\\
&&\ll \sum_{h}\frac{d^{k_1}(h)}{h^{3k_1/4}}
\sum_{s_1>(y/h)^{1/2}}\frac{d^2(s_1)}{s_1^{3/2}}
\sum_{s_{k_1}\gg (y/h)^{1/2}}\frac{d^2(s_{k_1})}{s_{k_1}^{3/2}}\\
&&\ll
\sum_{h}\frac{d^{k_1}(h)}{h^{3k_1/4}}(\frac{y}{h})^{-1/2+\varepsilon}
\ll y^{-1/2+\varepsilon}
\end{eqnarray*}
if we notice  $k_1\geq 3.$

If ${\cal G}_{k;l}={(I,J)},$ we have
\begin{equation}
U_2(d;y)\ll R_1^{(I,J)}(d;y) \ll y^{-1/2+\varepsilon}.
\end{equation}

If ${\cal G}_{k;l}\not={(I,J)},$
 we have
\begin{equation}
U_2(d;y)\ll R_1^{(I,J)}(d;y)+\sum_{\stackrel{(I_0,J_0)\in {\cal
G}_{k;l}} {(I_0,J_0)\not=
(I,J)}}R_1^{(I_0,J_0)}(d;y)R_2^{(I_0,J_0)}(d)\ll
y^{-1/2+\varepsilon}.
\end{equation}

Now Lemma 3.1 follows from (3.1)-(3.4).
\end{proof}

\section{\bf Proofs of Theorem 1 and Theorem 2}

Suppose $T\geq 10$ is a real number. It suffices for us to evaluate
the integral $\int_T^{2T}\Delta^k(x)dx.$ Suppose $y$ is a parameter
such that $T^\varepsilon<y\leq T^{1/3}.$ For any $T\leq x\leq 2T,$
define
\begin{eqnarray*}
&&{\cal R}_1={\cal R}_1(x,y): =(\sqrt 2\pi)^{-1}x^{1/4}\sum_{n\leq
y}\frac{d(n)}
{n^{3/4}}\cos(4\pi\sqrt{xn}-\frac{\pi}{4}),\\
&&{\cal R}_2={\cal R}_2(x,y): =\Delta(x)-{\cal R}_1.
\end{eqnarray*}
We shall show that the higher-power moment of ${\cal R}_2$ is small
and hence the  integral $\int_T^{2T}\Delta^k(x)dx$ can be well
approximated by $\int_T^{2T}{\cal R}_1^kdx,$ which is easy to
evaluate.

\subsection{\bf Evaluation of the integral $\int_T^{2T}{\cal R}_1^hdx$}\

Suppose $h\geq 3$ is any fixed integer. By the elementary formula
$$\cos a_1\cdots \cos a_h=
\frac{1}{2^{h-1}}\sum_{(i_1,\cdots,i_{h-1})\in \{0,1\}^{h-1}} \cos
(a_1+(-1)^{i_1}a_2+(-1)^{i_2}a_3+\cdots +(-1)^{i_{h-1}}a_h),$$ we
have
\begin{eqnarray*}
{\cal R}_1^h&&=(\sqrt 2\pi)^{-h}x^{\frac{h}{4}} \sum_{n_1\leq
y}\cdots \sum_{n_h\leq y}\frac{d(n_1)\cdots d(n_h)}{(n_1\cdots n_h
)^{3/4}}
\prod_{j=1}^h\cos(4\pi\sqrt{n_jx}-\pi/4)\\
&&=\frac{x^{\frac{h}{4}}}{(\sqrt 2\pi)^{h}2^{h-1}}
\sum_{(i_1,\cdots,i_{h-1})\in \{0,1\}^{h-1}} \sum_{n_1\leq y}\cdots
\sum_{n_h\leq y}
\frac{d(n_1)\cdots d(n_h)}{(n_1\cdots n_h)^{3/4}}\\
&&\hspace{10mm}\times \cos (4\pi\sqrt
x\alpha(n_1,\cdots,n_h;i_1,\cdots,i_{h-1})-\frac{\pi}{4}
\beta(i_1,\cdots,i_{h-1})),
\end{eqnarray*}
where
\begin{eqnarray*}
&&\alpha(n_1,\cdots,n_h;i_1,\cdots,i_{h-1})\\
&&\hspace{5mm}: =
\sqrt{n_1}+(-1)^{i_1}\sqrt{n_2}+(-1)^{i_2}\sqrt{n_3}+\cdots
+(-1)^{i_{h-1}}\sqrt
{n_h},\\
&&\beta(i_1,\cdots,i_{h-1}): =1+(-1)^{i_1}+(-1)^{i_2}+\cdots
+(-1)^{i_{h-1}}.
\end{eqnarray*}
Thus we can write
\begin{equation}
{\cal R}_1^h=\frac{1}{(\sqrt 2\pi)^{h}2^{h-1}}(S_1(x)+S_2(x)),
\end{equation}
where
\begin{eqnarray*}
S_1(x):&&= x^{h/4}\sum_{(i_1,\cdots,i_{h-1})\in
\{0,1\}^{h-1}}\cos(-\frac{\pi\beta}{4}) \sum_{\stackrel{n_j\leq
y,1\leq j\leq h}{\alpha= 0}}
\frac{d(n_1)\cdots d(n_h)}{(n_1\cdots n_h)^{3/4}},\\
S_2(x):&&= x^{h/4}\sum_{(i_1,\cdots,i_{h-1})\in \{0,1\}^{h-1}}
\sum_{\stackrel{n_j\leq y,1\leq j\leq h}{\alpha\not= 0}}
\frac{d(n_1)\cdots d(n_h)}{(n_1\cdots n_h)^{3/4}}
\cos(4\pi\alpha\sqrt x-\frac{\pi\beta}{4}),\\
\alpha: &&=\alpha(n_1,\cdots,n_h;i_1,\cdots,i_{h-1}),\\
\beta: &&=\beta(i_1,\cdots,i_{h-1}).\end{eqnarray*}

First consider the contribution of $S_1(x).$ We have
\begin{equation}
\int_T^{2T}S_1(x)dx= \sum_{(i_1,\cdots,i_{h-1})\in
\{0,1\}^{h-1}}\cos(-\frac{\pi\beta}{4}) \sum_{\stackrel{n_j\leq
y,1\leq j\leq h}{\alpha= 0}} \frac{d(n_1)\cdots d(n_h)}{(n_1\cdots
n_h)^{3/4}}\int_T^{2T}x^{\frac{h}{4}}dx.
\end{equation}

It is easily seen that if $\alpha=0,$ then $1\in
\{i_1,\cdots,i_{h-1}\}.$ Let $l=i_1+\cdots+i_{h-1},$ then we have
$$\sum_{\stackrel{n_j\leq y,1\leq j\leq h}{\alpha= 0}}
\frac{d(n_1)\cdots d(n_h)}{(n_1\cdots n_h)^{3/4}}=s_{h;l}(d;y),$$
where $s_{h;l}(d;y)$ was defined in last section. By Lemma 3.1 we
get
\begin{equation}
\int_T^{2T}S_1(x)dx=B_h^{*}(d)\int_T^{2T}x^{\frac{h}{4}}dx
+O(T^{1+h/4+\varepsilon}y^{-1/2}),
\end{equation}
where
$$B_h^{*}(d):=
\sum_{(i_1,\cdots,i_{h-1})\in
\{0,1\}^{h-1}}\cos(-\frac{\pi\beta}{4})
\sum_{\stackrel{(n_1,\cdots,n_h)\in {\Bbb N}^h}{\alpha= 0}}
\frac{d(n_1)\cdots d(n_h)}{(n_1\cdots n_h)^{3/4}}.$$

For any $(i_1,\cdots,i_{h-1})\in \{0,1\}^{h-1}\setminus
\{(0,\cdots,0)\},$ let
\begin{eqnarray*}
S(d;i_1,\cdots,i_{h-1}): &&=\sum_{\stackrel{(n_1,\cdots,n_h)\in
{\Bbb N}^h}{\alpha= 0}}
\frac{d(n_1)\cdots d(n_h)}{(n_1\cdots n_h)^{3/4}},\\
l(i_1,\cdots,i_{h-1}): &&=i_1+\cdots+i_{h-1}.
\end{eqnarray*}
It is easily seen that if
$l(i_1,\cdots,i_{h-1})=l(i_1^{\prime},\cdots, i_{h-1}^{\prime})$ or
$l(i_1,\cdots,i_{h-1})+l(i_1^{\prime},\cdots, i_{h-1}^{\prime})=h,$
then
$$S(d;i_1,\cdots,i_{h-1})=S(d;i_1^{\prime},\cdots,
i_{h-1}^{\prime})=s_{h;l(i_1,\cdots,i_{h-1})}(d).$$ From
$(-1)^i=1-2i(i=0,1)$ we also have
$$\beta(i_1,\cdots,i_{h-1})=h-2l(i_1,\cdots,i_{h-1}).$$
So we get
\begin{equation}
B_h^{*}(d)=\sum_{l=1}^{h-1}\sum_{l(i_1,\cdots,i_{h-1})=l}
\cos(-\frac{\pi\beta}{4})S(d;i_1,\cdots,i_{h-1})
\end{equation}
\begin{eqnarray*}
&&=\sum_{l=1}^{h-1}s_{h;l}(d)\cos
\frac{\pi(h-2l)}{4}\sum_{l(i_1,\cdots,i_{h-1})
=l}1\\
&&=\sum_{l=1}^{h-1}{h-1\choose l}s_{h;l}(d)\cos
\frac{\pi(h-2l)}{4}=B_h(d).
\end{eqnarray*}

Now we consider the contribution of $S_2(x).$ By Lemma 2.3 we get
\begin{equation}
\int_T^{2T}S_2(x)dx\ll T^{1/2+h/4} \sum_{(i_1,\cdots,i_{h-1})\in
\{0,1\}^{h-1}} \sum_{\stackrel{n_j\leq y,1\leq j\leq h}{\alpha\not=
0}} \frac{d(n_1)\cdots d(n_h)}{(n_1\cdots n_h)^{3/4}|\alpha|}.
\end{equation}

It suffices for us to estimate the sum
$$\Sigma(y;i_1,\cdots,i_{h-1})=\sum_{\stackrel{n_j\leq y,1\leq j\leq h}{\alpha\not= 0}}
 \frac{d(n_1)\cdots d(n_h)}{(n_1\cdots
n_h)^{3/4}|\alpha|}$$ for fixed $(i_1,\cdots,i_{h-1})\in
\{0,1\}^{h-1}.$

If $(i_1,\cdots,i_{h-1})=(0,\cdots,0),$ then
\begin{eqnarray*}
&&\Sigma(y;0,\cdots,0)\ll \sum_{n_j\leq y,1\leq j\leq h}
\frac{d(n_1)\cdots d(n_h)}{(n_1\cdots
n_h)^{3/4}(\sqrt{n_1}+\cdots+\sqrt{n_h} )}
\\
&&\ll \sum_{n_j\leq y,1\leq j\leq h}
\frac{d(n_1)\cdots d(n_h)}{(n_1\cdots n_h)^{3/4+1/2h}}\\
&&\ll y^{(h-2)/4}\log^h y,
\end{eqnarray*}
where we used the estimates
$$\sum_{n\leq u}\ll u\log u,\hspace{2mm} x_1+\cdots+x_h\gg (x_1\cdots x_k)^{1/h}
.$$

For  $(i_1,\cdots,i_{h-1})\not=(0,\cdots,0),$ by a splitting
argument we get that there exist a group of numbers $1<N_1,
N_2,\cdots, N_h<y$ such that
$$\Sigma(y;i_1,\cdots,i_{h-1})\ll \Sigma_1^{*}\log^h y,$$
where
$$\Sigma_1^{*}=\sum_{\stackrel{N_j<n_j\leq 2N_j,1\leq j\leq h}{\alpha\not= 0}}
\frac{d(n_1)\cdots d(n_h)}{(n_1\cdots n_h)^{3/4}|\alpha|}.$$ Without
loss of generality, we suppose $N_1\leq  N_2\leq \cdots\leq  N_h\leq
y.$

By Lemma 2.2 we have $|\alpha|\gg N_h^{-(2^{h-2}-2^{-1})}.$ Then by
a splitting argument and Lemma 2.4 we get for some
$N_h^{-(2^{h-2}-2^{-1})}\ll \Delta<y^{1/2}$ such that
\begin{eqnarray*}
\Sigma_1^{*}&&\ll \frac{y^\varepsilon}{(N_1\cdots N_h)^{3/4}\Delta}
{\cal A}(N_1,\cdots, N_h; i_1,\cdots,i_{h-1};\Delta)\\
&&\ll \frac{y^\varepsilon}{(N_1\cdots N_h)^{3/4}\Delta}
(\Delta N_h^{1/2}N_1\cdots N_{h-1}+N_1\cdots N_{h-1})\\
&&\ll y^\varepsilon (\frac{(N_1\cdots N_{h-1})^{1/4}}{N_h^{1/4}}
+\frac{(N_1\cdots N_{h-1})^{1/4}}{N_h^{3/4}\Delta})\\
&&\ll y^\varepsilon ( N_h^{(h-2)/4}+N_h^{b(h)})\ll
y^{b(h)+\varepsilon},
\end{eqnarray*}
where $b(h)$ was defined in Section 1.1. Thus we get
\begin{equation}
\int_T^{2T}S_2(x)dx\ll T^{1/2+h/4+\varepsilon}y^{b(h)}.
\end{equation}

Hence from (4.1)-(4.6) we get

{\bf Lemma 4.1.} For any fixed $h\geq 3,$ we have
\begin{equation}
\int_T^{2T}{\cal R}_1^hdx=\frac{B_h(d)}{(\sqrt
2\pi)^h2^{h-1}}\int_T^{2T} x^{h/4}dx
+O(T^{1+h/4+\varepsilon}y^{-1/2}+T^{1/2+h/4+\varepsilon}y^{b(h)}).
\end{equation}

\subsection{\bf Higher-power moments of ${\cal R}_2$}\

We first study the mean-square of ${\cal R}_2.$ We begin with the
truncated
 Voronoi's formula[9, equation (2.25)]
\begin{equation}
\Delta(x)=(\pi\sqrt 2)^{-1}x^{\frac{1}{4}}\sum_{n\leq
N}\frac{d(n)}{n^{3/4}} \cos(4\pi\sqrt{nx}-\pi/4)
+O(x^{1/2+\varepsilon}N^{-1/2}),
\end{equation}
where $1< N\ll x.$ Taking $N=T,$ we get
\begin{eqnarray*}
{\cal R}_2 &&=(\pi\sqrt 2)^{-1}x^{\frac{1}{4}}\sum_{y<n\leq
T}\frac{d(n)}{n^{3/4}} \cos(4\pi\sqrt{nx}-\pi/4)
+O(T^{\varepsilon})\\
&&\ll |x^{\frac{1}{4}}\sum_{y<n\leq T}\frac{d(n)}{n^{3/4}}
e(2\sqrt{nx})|+T^\varepsilon,
\end{eqnarray*}
which implies
\begin{equation}
\int_T^{2T}{\cal R}_2^2dx \ll
T^{1+\varepsilon}+\int_T^{2T}|x^{\frac{1}{4}}\sum_{y<n\leq
T}\frac{d(n)}{n^{ 3/4}} e(2\sqrt{nx})|^2dt
\end{equation}
\begin{eqnarray*}
&&\ll T^{1+\varepsilon}+T^{3/2}\sum_{y<n\leq
T}\frac{d^2(n)}{n^{3/2}}
+T\sum_{y<m<n\leq T}\frac{d(n)d(m)}{(mn)^{3/4}(\sqrt n-\sqrt m)}\\
&&\ll T^{1+\varepsilon}+\frac{T^{3/2}\log^3 T}{y^{1/2}} \ll
\frac{T^{3/2}\log^3 T}{y^{1/2}},
\end{eqnarray*}
where we used the estimates
\begin{eqnarray*}
&&\sum_{n\leq u}d^2(n)\ll u\log^3 u,\\
&&\sum_{y<m<n\leq T}\frac{d(n)d(m)}{(mn)^{3/4}(\sqrt n-\sqrt m)}\ll
T^\varepsilon.
\end{eqnarray*}

Now suppose $y$ satisfies $y^{2b(K_0)}\leq T.$ Hence from Lemma 4.1
we get that
$$ \int_T^{2T}|{\cal R}_1|^{K_0}dx\ll T^{1+K_0/4+\varepsilon},$$
which implies
\begin{equation}
 \int_T^{2T}|{\cal R}_1|^{A_0}dx\ll T^{1+A_0/4+\varepsilon}
\end{equation}
since $A_0\leq K_0.$ From (1.10) and (4.10) we get
\begin{equation}
 \int_T^{2T}|{\cal R}_2|^{A_0}dx\ll
 \int_T^{2T}(|\Delta(x)|^{A_0}+|{\cal R}_1|^{A_0})dx\ll T^{1+A_0/4+\varepsilon}.
\end{equation}

For any $2<A<A_0,$ from (4.9), (4.11) and  H\"older's inequality we
get that
\begin{equation}
 \int_T^{2T}|{\cal R}_2|^{A}dx=
\int_T^{2T}|{\cal
R}_2|^{\frac{2(A_0-A)}{A_0-2}+\frac{A_0(A-2)}{A_0-2}}dx
\end{equation}
$$\ll \left( \int_T^{2T}{\cal R}_2^2dx\right)^{\frac{A_0-A}{A_0-2}}
 \left( \int_T^{2T}|{\cal R}_2|^{A_0}dx\right)^{\frac{A-2}{A_0-2}}
\ll T^{1+\frac{A}{4}+\varepsilon}y^{-\frac{A_0-A}{2(A_0-2)}}.$$

Namely, we have the following Lemma 4.2.

{\bf Lemma 4.2.} Suppose $T^\varepsilon\leq y\leq T^{1/2b(K_0)},$
$2<A<A_0,$  then
\begin{equation}
 \int_T^{2T}|{\cal R}_2|^{A}dx\ll
T^{1+\frac{A}{4}+\varepsilon}y^{-\frac{A_0-A}{2(A_0-2)}}.
\end{equation}

\subsection{\bf Proof of Theorem 1.}\

Suppose $3\leq k\leq K(A_0)$ and  $T^\varepsilon\leq y\leq
T^{1/2b(K_0)}.$ By the elementary formula $(a+b)^k-a^k\ll
|a^{k-1}b|+|b|^k,$ we get
\begin{equation}
\int_T^{2T}\Delta^k(x)dx=\int_T^{2T}{\cal R}_1^kdx
+O(\int_T^{2T}|{\cal R}_1^{k-1}{\cal R}_2|dx) +O(\int_T^{2T}|{\cal
R}_2|^kdx).
\end{equation}

If $k-1<A_0/2,$ then from (4.9), (4.10) and Cauchy's inequality we
get
$$
\int_T^{2T}|{\cal R}_1^{k-1}{\cal R}_2|dx \ll
\left(\int_T^{2T}|{\cal R}_1|^{2(k-1)}dx\right)^{1/2}
\left(\int_T^{2T}|{\cal R}_2|^2dx\right)^{1/2} \ll
T^{1+k/4+\varepsilon}y^{-1/4}.$$

If $k-1\geq A_0/2,$ then from (4.10), Lemma 4.2 and H\"older's
inequality we get
\begin{eqnarray*}
\int_T^{2T}|{\cal R}_1^{k-1}{\cal R}_2|dx &&\ll
\left(\int_T^{2T}|{\cal R}_1|^{A_0}dx\right)^{\frac{k-1}{A_0}}
\left(\int_T^{2T}|{\cal R}_2|^{\frac{A_0}{A_0-k+1}}dx\right)^
{\frac{A_0-k+1}{A_0}}\\
&& \ll T^{1+k/4+\varepsilon}y^{-(A_0-k)/2(A_0-2)}.\end{eqnarray*}
Thus we have
\begin{equation}
\int_T^{2T}|{\cal R}_1^{k-1}{\cal R}_2|dx +\int_T^{2T}|{\cal
R}_2|^kdx\ll T^{1+k/4+\varepsilon}y^{-\sigma(k,A_0)},
\end{equation}
where $\sigma(k,A_0)$ was defined in Section 1.1.

From (4.14) and (4.15) we get
\begin{equation}
\int_T^{2T}\Delta^k(x)dx=\int_T^{2T}{\cal R}_1^kdx
+O(T^{1+k/4+\varepsilon}y^{-\sigma(k,A_0)}).
\end{equation}

Now take $y=T^{1/2b(K_0)}.$ From Lemma 4.1 and (4.16) we get
\begin{equation}
\int_T^{2T}\Delta^k(x)dx=\frac{B_k(d)}{(\sqrt
2\pi)^k2^{k-1}}\int_T^{2T}x^{k/4}d x
+O(T^{1+k/4-\sigma(k,A_0)/2b(K_0)+\varepsilon})
\end{equation}
$$=\frac{B_k(d)}{(\sqrt 2\pi)^k2^{k-1}}\int_T^{2T}x^{k/4}dx
+O(T^{1+k/4-\delta_1(k,A_0)+\varepsilon}).$$ Theorem 1 follows from
(4.17) immediately.

\subsection{\bf Proof of Theorem 2}\

Suppose $T^\varepsilon\leq y\leq T^{1/3}.$ By the truncated
Voronoi's formula (4.8), we have
$${\cal R}_2=(\sqrt 2\pi)^{-1}x^{1/4}\sum_{y<n\leq N}\frac{d(n)}{n^{3/4}}
\cos(4\pi\sqrt{nx}-\pi/4)+O(x^{1/2+\varepsilon}N^{-1/2}),$$
 where $y<N\ll T.$
Using Ivi\'c's large-value technique directly to ${\cal R}_2$
without modifications , we  get that the estimate
\begin{equation}
 \int_T^{2T}|{\cal R}_2|^{A_0}dx\ll T^{1+A_0/4+\varepsilon}
\end{equation}
 holds with $A_0=184/19,   T^\varepsilon\leq y\leq T^{1/3}.$
We omit the details since it is  completely the same as that of
Ivi\'c. Combining (1.10) we get that
\begin{equation}
 \int_T^{2T}|{\cal R}_1|^{A_0}dx\ll T^{1+A_0/4+\varepsilon}
\end{equation}
 holds with $A_0=184/19, T^\varepsilon\leq y\leq T^{1/3}.$

By the same argument as in last subsection , we get that for
$T^\varepsilon\leq y\leq T^{1/3},$
\begin{equation}
\int_T^{2T}\Delta^k(x)dx=\int_T^{2T}{\cal R}_1^kdx
+O(T^{1+k/4+\varepsilon}y^{-\sigma(k,184/19)}).
\end{equation}
 Take $y=T^{1/(2b(k)+2\sigma(k,184/19))}.$
From Lemma 4.1 again we get
\begin{equation}
\int_T^{2T}\Delta^k(x)dx=\frac{B_k(d)}{(\sqrt
2\pi)^k2^{k-1}}\int_T^{2T}x^{k/4}d x+
O(T^{1+k/4-\frac{\sigma(k,184/19)}{
2b(k)+2\sigma(k,184/19)}+\varepsilon})
\end{equation}
$$=\frac{B_k(d)}{(\sqrt 2\pi)^k2^{k-1}}\int_T^{2T}x^{k/4}dx+
O(T^{1+k/4-\delta_2(k,184/19)+\varepsilon}).$$ And Theorem 2
follows.

\section{\bf Proofs of other Theorems}

$P(x)$ has the following truncated Voronoi's formula
\begin{equation}
P(x)=-\frac{1}{\pi}\sum_{n\leq
N}r(n)n^{-3/4}x^{1/4}\cos(4\pi\sqrt{nx}+\pi/4)
+O(x^{1/2+\varepsilon}N^{-1/2})
\end{equation}
for $1\leq N\ll x,$ which follows from Lemma 3 of M\"uller[16].
$A(x)$ has the following truncated Voronoi's formula
\begin{equation}
A(x)=\frac{1}{\pi\sqrt 2}x^{\kappa/2-1/4}\sum_{n\leq
N}a(n)n^{-\kappa/2-1/4} \cos(4\pi\sqrt{nx}-\pi/4)
\end{equation}
$$+O(x^{\kappa/2+\varepsilon}N^{-1/2})$$
for $1\leq N\ll x,$ which is a special case of Theorem 1.1 of
Jutila[13]. So by the same way as in last section, we get Theorem 3
and Theorem 4 .

Now we prove Theorem 5. We shall follow Ivi\'c[10]. Define
$$\Delta^{*}(x): =\frac{1}{2}\sum_{n\leq 4x}(-1)^nd(n)-x(\log x+2\gamma-1), x>0.
$$

Jutila[12] proved that
\begin{equation}
\int_0^T(E(t)-2\pi\Delta^{*}(\frac{t}{2\pi}))^2dt\ll T^{4/3}\log^3
T,
\end{equation}
which means that $E(t)$ is well approximated by
$2\pi\Delta^{*}(\frac{t}{2\pi})$ at least in the mean square sense.

Suppose $A_0>9$ is a real number such that  both of (1.10) and (1.32
) hold. Since (see Jutila[11])
$$\Delta^{*}(x)=-\Delta(x)+2\Delta(2x)-\frac{1}{2}\Delta(4x),$$
from (1.10) we get
\begin{equation}
\int_0^T|\Delta^{*}(t)|^{A_0}dt\ll T^{1+A_0/4+\varepsilon}.
\end{equation}
Then from (1.32), (5.3), (5.4) and H\"older's inequality we get for
any $3\leq k<A_0$ that
\begin{equation}
\int_0^TE^k(t)dt-(2\pi)^{k+1}\int_0^{\frac{T}{2\pi}}(\Delta^{*}(t))^kdt
\end{equation}
\begin{eqnarray*}
&&=\int_0^T\left(E^k(t)-(2\pi\Delta^{*}(\frac{t}{2\pi}))^k\right)dt\\
&&\ll \int_0^T|E(t)-2\pi\Delta^{*}(\frac{t}{2\pi})|
\left(|E(t)|^{k-1}+|\Delta^{*}(\frac{t}{2\pi})|^{k-1}\right)dt\\
&&\ll T^{1+k/4-\sigma(k,A_0)/3+\varepsilon},
\end{eqnarray*}
where $\sigma(k,A_0)$ was defined in Section 1.1.

From (5.5) the problem is reduced to evaluating the integral
$\int_0^T(\Delta^{*}(t))^kdt.$
 For $1\ll N\ll x,$ we have[10, equation (7)]
\begin{equation}
\Delta^{*}(x)=\frac{1}{\pi\sqrt 2}\sum_{n\leq
N}(-1)^nd(n)n^{-3/4}x^{1/4} \cos(4\pi\sqrt{nx}-\pi/4)
+O(x^{1/2+\varepsilon}N^{-1/2}),
\end{equation}
which is similar to (4.8). Let $d^{*}(n)=(-1)^nd(n).$ Then  by the
same way as i n the proof of Theorem 1, we get that the asymptotic
formula
\begin{equation}
\int_1^T(\Delta^{*}(t))^kdt=\frac{B_k(d^{*})}{(1+k/4)2^{3k/2-1}\pi^k}T^{1+k/4}
+O(T^{1+k/4-\delta_1(k,A_0)+\varepsilon})
\end{equation}
holds for any  $3\leq k<A_0.$

We shall use the following Lemma 5.1.

{\bf Lemma 5.1} Suppose $1\leq l<k$ are fixed integers,
$(n_1,\cdots, n_k)\in {\Bbb N}^k.$ If
$$\sqrt{n_1}+\cdots+\sqrt{n_l}=\sqrt{n_{l+1}}+\cdots+\sqrt{n_k}$$
holds , then $2|(n_1+\cdots+n_k).$

\begin{proof}
For any $n\in {\Bbb N},$ let $h(n)$ denote the squarefree part of
$n.$ Let ${\cal S}=\{h(n_1),$ $\cdots,$ $ h(n_k)\}$ $\bigcap{\Bbb
N}$ and $s=\#{\cal S}.$ For convenience, write
$${\cal S}=\{h_1,\cdots, h_s\}, I=\{1,\cdots, l\}, J=\{l+1,\cdots,k\}.$$
From Lemma 2.1 we can write $I=\bigcup_{e=1}^sI_e,
J=\bigcup_{e=1}^sJ_e$ such that for each $1\leq e\leq s,$ we have
$$\sum_{i\in I_e}\sqrt{n_i}=\sum_{j\in J_e}\sqrt{n_j}$$
and that all $n_i(i\in I_e)$ and $n_j(j\in J_e)$ have the same
squarefree part $ h_e.$ Namely we have ($1\leq e\leq s$)
$$n_i=m_i^2h_e(i\in I_e),n_j=m_j^2h_e(i\in I_e),
\sum_{i\in I_e}m_i=\sum_{j\in J_e}m_j .$$
 Thus we get
\begin{eqnarray*}
&&n_1+\cdots+n_k=\sum_{e=1}^s
(\sum_{i\in I_e}n_i+\sum_{j\in J_e}n_j)\\
&&=\sum_{e=1}^s(\sum_{i\in I_e}m_i^2h_e+\sum_{j\in J_e}m_j^2h_e)
\equiv\sum_{e=1}^s(\sum_{i\in I_e}m_i+\sum_{j\in J_e}m_j)h_e\\
&&=2\sum_{e=1}^sh_e\sum_{i\in I_e}m_i\equiv 0(mod\hspace{2mm}2),
\end{eqnarray*}
where we used the simple congruence $n^2\equiv n(mod\hspace{2mm}2).$
\end{proof}

From Lemma 5.1 we get for any $1\leq l<k$ that
\begin{eqnarray*}
s_{k;l}(d^{*})&& =\sum_{\sqrt{n_1}+\cdots
+\sqrt{n_l}=\sqrt{n_{l+1}}+\cdots +\sqrt{n_k}}
(-1)^{n_1+\cdots+n_k}\frac{d(n_1)\cdots d(n_k)}{(n_1\cdots n_k)^{3/4}}\\
&&=\sum_{\sqrt{n_1}+\cdots +\sqrt{n_l}=\sqrt{n_{l+1}}+\cdots
+\sqrt{n_k}}
\frac{d(n_1)\cdots d(n_k)}{(n_1\cdots n_k)^{3/4}}\\
&&=s_{k;l}(d).
\end{eqnarray*}
Whence we get
\begin{equation}
B_k(d^{*})=B_k(d).
\end{equation}

From (5.5), (5.7) and (5.8) we get (1.33).

Similar to Theorem 2, we can prove that the asymptotic formula
\begin{equation}
\int_1^T(\Delta^{*}(t))^kdt=\frac{B_k(d)}{(1+k/4)2^{3k/2-1}\pi^k}T^{1+k/4}
+O(T^{1+k/4-\delta_2(k,576/61)+\varepsilon})
\end{equation}
holds for any  $3\leq k\leq 9,$ which combined with (5.5)  yields
the second part of Theorem 3.

{\bf Note added in proof}: Recently M. N. Huxley (Exponential sums
and Lattice points III, Proc. London Math. Soc., Vol.{\bf
87}(3)(2003),
 591-609) proved $$\Delta(x)\ll x^{131/416}(\log x)^{26947/8320},$$ which implies that
the formula (1.10) holds for $A_0=262/27.$ The exponent
$\delta_2(k,184/19) $ in Theorem 2 then can be improved to
$\delta_2(k,262/27)$ for $k=6,7,8,9.$ The author deeply thanks
Professor A. Schinzel  for informing me M. N. Huxley's new result.


\medskip



\begin{thebibliography}{99}
 \bibitem{s1}F. V. Atkinson, the mean value of the Riemann zeta-function,  Acta
Math.{\bf 81}(1949), 353-376.

\bibitem{s2}A. S. Besicovitch, On the linear independence of fractional powers
of integers, J. London Math. Soc.{\bf 15}(1940),3-6.

\bibitem{s2}Cai Yingchun, On the third and fourth power moments of Fourier
coefficients of cusp forms.  Acta Math. Sinica (N.S.) 13 (1997), no.
4, 443--452 .




\bibitem{s3} P. Deligne, La conjecture de Weil I. Inst. Hautes Etudes Sci.
Publ. Math.{\bf 43}(1974), 273-307.


\bibitem{s10}D. R. Heath-Brown, The distribution and moments of the error term i
n the Dirichlet divisor problem, Acta Arith. (1992), 389-415.

\bibitem{s7}M. N. Huxley,
Area, lattice points, and exponential sums. London Math. Soc.
Monographs.
 New Series, 13. Oxford University Press, 1996.

\bibitem{s1}A. Ivi\'c, The Riemann zeta-function. John. Wiley and Sons, 1985.


\bibitem{s2}A. Ivi\'c, Lectures on mean values of the Riemann zeta-function, Lec
tures On Math. and Physics 82, Tata Inst. Fund. Res., Bombay, 1991.

\bibitem{s1}  A. Ivi\'c, Large values of certain number-theoretic error terms
Acta Arith. {\bf 56}(1990), 135-159.

\bibitem{s4}A. Ivi\'c, On some problems involving the mean square of
 $\zeta(\frac 12+it)$. Bull. Cl. Sci. Math. Nat. Sci. Math. No. 23 (1998), 71--7
6.


\bibitem{s2}M. Jutila, Riemann's zeta-function and the divisor problem
Ark. Mat. {\bf 21}(1983), 75-96 and II,ibid.
 {\bf 31}(1993), 61-70.

\bibitem{s1}M. Jutila, On a formula of Atkinson.
Topics in classical number theory, Vol. I, II (Budapest, 1981),
807--823, Colloq. Math. Soc. Janos Bolyai, 34, North-Holland,
Amsterdam,
 1984.

\bibitem{s2}M. Jutila, Lectures on a method in the theory of exponential sums.
Tata Inst. of Fund. Res. Lectures on Math. and Physics, 80.
 Bombay, 1987.


\bibitem{s2}I. Katai, The number of lattice points in a circle. (Russian)
Ann. Univ. Sci. Budapest. E\"otv\"os Sect. Math. {\bf 8}(1965), 39--
60.



\bibitem{s1}T. Meurman, On the mean square of the Riemann zeta-function. Quart.
J. Math. Oxford Ser. (2) 38 (1987), no. 151, 337--343.

\bibitem{s1}W. M\"uller, On the asymptotic behaviour of the ideal counting
function in quadratic number fields. Monatsh. Math. 108 (1989), no.
4, 301--323.


\bibitem{s3}K. C. Tong, On divisor problem III, Acta math. Sinica {\bf 6}
(1956), 515-541.

\bibitem{s4}Kai-Man Tsang, Higher-power moments of $\Delta(x), E(t)$ and $P(x)$,
Proc. London Math. Soc.(3){\bf 65}(1992), 65-84.

\bibitem{s4}G. Voronoi, Surune fonction transcendante et ses applications a la
sommation de quelques s\'eries, Ann. Sci. \'Ecole Norm Sup. (3){\bf
21}(1904), 207-267, 459-533.

\bibitem{s3} Wenguang Zhai, On higher-power moments of $\Delta(x),$  Acta Arith. Vol.{\bf 112}(2004), 1-24.

\end{thebibliography}
\end{document}